\documentclass[fleqn]{mat01}
\usepackage{times,mathtimy,amssymb,latexsym}
\begin{document}

\setcounter{page}{185} \firstpage{185}

\newtheorem{theore}{Theorem}
\renewcommand\thetheore{\arabic{section}.\arabic{theore}}
\newtheorem{theor}[theore]{\bf Theorem}
\newtheorem{lem}[theore]{Lemma}

\renewcommand\theequation{\arabic{section}.\arabic{equation}}

\title{Weighted composition operators on weighted Bergman spaces of
bounded symmetric domains}

\markboth{Sanjay Kumar and Kanwar Jatinder Singh}{Bergman spaces
of bounded symmetric domains}

\author{SANJAY KUMAR and KANWAR JATINDER SINGH}

\address{Department of Mathematics, University of Jammu, Jammu~180 006, India\\
\noindent E-mail: sks$_{-}$jam@yahoo.co.in;
kunwar752000@yahoo.co.in}

\volume{117}

\mon{May}

\parts{2}

\pubyear{2007}

\Date{MS received 7 December 2005}

\begin{abstract}
In this paper, we study the weighted compositon operators on
weighted Bergman spaces of bounded symmetric domains. The
necessary and sufficient conditions for a weighted composition
operator $W_{\varphi,\psi}$ to be bounded and compact are studied
by using the Carleson measure techniques. In the last section, we
study the Schatten $p$-class weighted composition operators.
\end{abstract}


\maketitle

\section{Introduction}

Let $\Omega$ be a bounded symmetric domain in $\mathbf{C}^{n}$
with Bergman kernel $K(z,\omega).$ We assume that $\Omega$ is in
its standard representation and the volume measure ${\rm d}V$ of
$\Omega$ is normalised so that $K(z,0) = K(0,\omega),$ for all $z$
and $\omega$ in $\Omega.$ By Theorem~5.7 of \cite{koranyi} and
using the polar coordinates representation, there exists a
positive number $\varepsilon_{\Omega}$ such that for $\lambda <
\varepsilon_{\Omega},$ we have
\begin{equation*}
C_{\lambda} = \int_{\Omega} K(z,z)^{\lambda} {\rm d}V(z) < \infty.
\end{equation*}
For each $\lambda < \varepsilon_{\Omega},$ define ${\rm
d}V_{\lambda}(z) = C_{\lambda}^{-1} K(z,z)^{\lambda} {\rm d}V(z).
$ Then $\{{\rm d}V_{\lambda}\}$ defines a weighted family of
probability measures on $\Omega.$ Also, throughout the paper
$\lambda < \varepsilon_{\Omega}$ is fixed. We define the weighted
Bergman spaces $A^{p}_{\lambda}(\Omega, {\rm d}V_{\lambda}),$ on
$\Omega$, as the set of all holomorphic functions $f$ on $\Omega$
so that
\begin{equation*}
\| f\|_{A^{p}_{\lambda}} = \left(\int_{\Omega} |f|^{p} {\rm
d}V_{\lambda} \right)^{\frac{1}{p}} < \infty.
\end{equation*}
Note that $A^{p}_{\lambda}(\Omega, {\rm d}V_{\lambda})$ is a
closed subspace of $L^{p}(\Omega, {\rm d}V_{\lambda}). $ For $
\lambda = 0, A^{p}_{\lambda}(\Omega, {\rm d}V_{\lambda}) $ is just
the usual Bergman space. For $p = 2,$ there is an orthogonal
projection $P_{\lambda}$ from $L^{2}(\Omega, {\rm d}V_{\lambda})$
onto $A^{2}_{\lambda}(\Omega, {\rm d}V_{\lambda}) $ given by
\begin{equation*}
P_{\lambda}f(z) = \int_{\Omega} K_{\lambda}(z,\omega) f(\omega)
{\rm d}V_{\lambda}(\omega),
\end{equation*}
where $K_{\lambda}(z,\omega)= K(z,\omega)^{1 - \lambda} $ is the
reproducing kernel for $A^{2}_{\lambda}(\Omega, {\rm
d}V_{\lambda}). $

Suppose, $\varphi, \psi$ are holomorphic mappings defined on $
\Omega$ such that $\varphi(\Omega) \subseteq \Omega. $ Then the
weighted composition operator $W_{\varphi,\psi}$ is defined as
\begin{equation*}
W_{\varphi,\psi}f(z) = \psi(z)f(\varphi (z)), \;{\rm for \; all
\;} f {\rm \; holomorphic \; in} \; \Omega, \;{\rm and} \; z \in
\Omega.
\end{equation*}
For the study of weighted composition operators one can refer to
\cite{con,cuc,kumar,smith} and references therein. Recently, Smith
\cite{wsmith} has made a nice connection between the Brennan's
conjecture and weighted composition operators. He has shown that
Brennan's conjecture is equivalent to the existence of self-maps
of unit disk that make certain weighted composition operators
compact.

We know that the bounded symmetric domain $\Omega$ in its standard
representation with normalised volume measure, the kernel function
$K(\cdot,\cdot)$ has the following special\break
properties:\vspace{-.3pc}
\begin{enumerate}
\renewcommand\labelenumi{(\arabic{enumi})}
\leftskip .15pc
\item $K(0,a) = 1 = K(a,0). $
\item $K(z,a) \neq 0$ (for all $z$ in $\Omega$ and $a$ in $\overline{\Omega}$).
\item $\lim_{a \rightarrow \partial \Omega} K(a,a) = \infty. $
\item $K(z,a)^{-1}$ is a smooth function on $\mathbf{C}^{n}
\times \mathbf{C}^{n}.$ Here $\overline{\Omega} $ denotes the
closure of $\Omega$ in $\mathbf{C}^{n}$ and $\partial \Omega$ is
the topological boundary, see \cite{hua} and \cite{koecher} for
details.\vspace{-.5pc}
\end{enumerate}

For any $a$ in $\Omega$, let
\begin{equation*}
k_{a}(z) = \frac{K(z,a)}{\sqrt{K(a,a)}}.
\end{equation*}
The $k_{a}'s$ are called normalised reproducing kernels for $
A^{2}(\Omega, {\rm d}V).$ They are unit vectors in $A^{2}(\Omega,
{\rm d}V).$ Moreover, $k_{a}^{1 - \lambda}$ is a unit vector of
$A^{2}_{\lambda}(\Omega, {\rm d}V_{\lambda}) $ for any $a$ in
$\Omega. $

Let $\mu$ be a finite complex Borel measure on $\Omega. $ Then the
Berezin transform of measure $\mu,$ denoted by $
\tilde{\mu}_{\lambda},$ is defined as
\begin{equation*}
\tilde{\mu}_{\lambda}(z) = \int_{\Omega} |k_{z}(\omega) |^{2 (1 -
\lambda)} {\rm d}\mu(\omega), \;\; z \in \Omega.
\end{equation*}
We will denote by $\beta(z,\omega)$ the Bergman distance function
on $\Omega.$ For $z \in \Omega$ and $r
> 0,$ let
\begin{equation*}
E(z,r) = \{\omega \in \Omega{\rm :}\ \beta(z,\omega)< r\}.
\end{equation*}
We denote by $|E(z,r)|$ the normalised volume of $E(z,r), $ that
is,
\begin{equation*}
|E(z,r)| = \int_{E(z,r)} {\rm d}V(\omega).
\end{equation*}
Given a finite complex Borel measure $\mu$ on $\Omega, $ we define
a function $\hat{\mu}_{r}$ on $\Omega $ by
\begin{equation*}
\hat{\mu}_{r}(z) = \frac{\mu(E(z,r))}{|E(z,r)|^{1-\lambda}}, \quad
z \in \Omega.
\end{equation*}

\begin{theor}[{\cite{zhu1}}]
Take $1 \le p < \infty.$ Let $\mu$ be a finite positive Borel
measure on $\Omega.$ Then the following conditions are
equivalent{\rm :}\vspace{-.3pc}
\begin{enumerate}
\renewcommand\labelenumi{\rm (\arabic{enumi})}
\leftskip .15pc
\item The inclusion map $i\hbox{\rm :}\ A^{p}_{\lambda}(\Omega,
{\rm d}V_{\lambda}) \rightarrow L^{p}(\Omega, {\rm d}V_{\lambda})
$ is bounded.
\item The Berezin transform $\tilde{\mu}_{\lambda}$ is bounded.
\item $\hat{\mu}_{r}$ is bounded on $\Omega$ for all {\rm (}or some{\rm
)} $r > 0. $
\item $\{\hat{\mu}_{r}(a_{n})\}$ is a bounded sequence{\rm ,}
where $\{a_{n}\}$ is some bounded sequence in $\Omega $
independent of $\mu. $\vspace{-.5pc}
\end{enumerate}
\end{theor}
A positive Borel measure $\mu$ satisfying any one of the
conditions of Theorem~1.1 is called a Carleson measure on the
weighted Bergman space $A^{p}_{\lambda}(\Omega, {\rm
d}V_{\lambda}) $

\begin{theor}[{\cite{zhu1}}] Take $1 \le p < \infty.$ Let $\mu$ be
a finite positive Borel measure on $\Omega.$ Then the following
conditions are equivalent{\rm :}\vspace{-.3pc}
\begin{enumerate}
\renewcommand\labelenumi{\rm (\arabic{enumi})}
\leftskip .15pc
\item The inclusion map $i\hbox{\rm :}\ A^{p}_{\lambda}(\Omega,
{\rm d}V_{\lambda}) \rightarrow L^{p}(\Omega, {\rm d}V_{\lambda})$
is compact.

\item The Berezin transform $\tilde{\mu}_{\lambda}(z) \rightarrow
0$ as $z \rightarrow \partial \Omega.$

\item $\hat{\mu}_{r}(z) \rightarrow 0$ as $z \rightarrow
\partial \Omega$ for all {\rm (}or some{\rm )} $r > 0.$

\item $\hat{\mu}_{r}(a_{n}) \rightarrow 0$ as $n \rightarrow
\infty. $\vspace{-.5pc}
\end{enumerate}
\end{theor}
A positive Borel measure $\mu$ satisfying any one of the
conditions of Theorem~1.2 is called a vanishing Carleson measure
on the weighted Bergman space $A^{p}_{\lambda}(\Omega, {\rm
d}V_{\lambda}). $

A positive compact operator $T$ on $A^{2}_{\lambda}(\Omega, {\rm
d}V_{\lambda})$ is in the trace class if
\begin{equation*}
{\rm tr}(T) = \sum_{n = 1}^{\infty} \langle Te_{n}, e_{n} \rangle
< \infty,
\end{equation*}
for some (or all) orthonormal basis $\{e_{n}\}$ of
$A^{2}_{\lambda}(\Omega, {\rm d}V_{\lambda}). $

Take $1 \le p < \infty,$ and $T$ be a compact linear operator on
$A^{2}_{\lambda}(\Omega, {\rm d}V_{\lambda}).$ Then we say that
$T$ belongs to the Schatten $p$-class $ S_{p} $ if $(T^{*}
T)^{{p}/{2}}$ is in the trace class. Also, the $S_{p} $ norm of
$T$ is given as
\begin{equation*}
\|T\|_{S_{p}} = [{\rm tr}((T^{*} T)^{{p}/{2}})]^{{1}/{p}}.
\end{equation*}
Moreover, $S_{p}$ is a two-sided ideal of the full algebra
$B(A^{2}_{\lambda}(\Omega, {\rm d}V_{\lambda}))$ of bounded linear
operators on $A^{2}_{\lambda}(\Omega, {\rm d}V_{\lambda}). $

\section{Preliminaries}

\setcounter{theore}{0}

\setcounter{equation}{0}

To make the paper self-contained we state the following lemmas.

\begin{lem}\hskip -.3pc{\rm \cite{kranz}.}\ \
The Bergman kernel $K(z,\omega)$ is conjugate symmetric. That
is{\rm ,}
\begin{equation*}
K(z,\omega) = \overline{K(\omega,z)};
\end{equation*}
and it is the function of $z$ in $ A^{2}_{\lambda}(\Omega, {\rm
d}V_{\lambda}). $
\end{lem}

\begin{lem}\hskip -.3pc{\rm \cite{zhu1}.}\ \
Let $T$ be a trace class operator or a positive operator on $
A^{2}_{\lambda}(\Omega, {\rm d}V_{\lambda}).$ Then
\begin{equation*}
{\rm tr}(T) = \int_{\Omega} \langle T K_{z}^{1-\lambda},
K_{z}^{1-\lambda} \rangle_{A^{2}_{\lambda}} {\rm d}\sigma(z),
\end{equation*}
where ${\rm d}\sigma(z) = K(z,z) {\rm d}V(z). $
\end{lem}

\begin{lem}\hskip -.3pc{\rm \cite{jafari}.}\ \
Let $\Omega$ be a bounded symmetric domain in $\mathbf{C}^{n}$ and
let $W_{\varphi,\psi}$ be a weighted composition operator on
$A^{2}_{\lambda}(\Omega, {\rm d}V_{\lambda}).$ Then
\begin{equation*}
\|W_{\varphi,\psi}^{*}\|_{A^{2}_{\lambda}} =
\;\;|\psi(z)|\left[\frac{K(\varphi(z),\varphi(z))}{K(z,z)}
\right]^{\frac{(1 - \lambda)}{2}}.
\end{equation*}
\end{lem}

\begin{lem}\hskip -.3pc{\rm \cite{arazy}.}\ \
Suppose $T$ is a positive operator on a Hilbert space $ H $ and
$x$ is a unit vector in $H$. Then\vspace{-.3pc}
\begin{enumerate}
\renewcommand\labelenumi{\rm (\arabic{enumi})}
\leftskip .15pc
\item $\langle T^{p}x, x \rangle \ge \langle Tx, x \rangle^{p} \;\;
{\rm for \;\; all}\;\; p \ge 1. $

\item $\langle T^{p}x, x \rangle \le \langle Tx, x \rangle^{p} \;\;
{\rm for \;\; all}\;\; 0 < p \le 1. $\vspace{-.5pc}
\end{enumerate}
\end{lem}

\begin{lem}\hskip -.3pc{\rm ({\it Theorem}~18.11({\it f}),
{\it p}.~89 {\it of} \cite{conway}).}\ \ Suppose $T$ is a bounded
linear operator on a Hilbert space $H$. If $T \in S_{p}, $ then
\begin{equation*}
{\rm tr}(| T |^{p}) = {\rm tr}(| T^{*} |^{p}).
\end{equation*}
\end{lem}

\begin{lem}\hskip -.3pc{\rm \cite{zhu1}.}\ \
Take $1 \le p < \infty.$ Suppose $\mu \ge 0$ is a finite Borel
measure on $\Omega.$ Then $\mu$ is a Carleson measure on
$A^{p}_{\lambda}(\Omega, {\rm d}V_{\lambda})$ if and only if
$\mu(E(z,r))/| E(z,r) |^{1 - \lambda}$ is bounded on $ \Omega $
{\rm (}as a function of $z${\rm )} for all {\rm (}or some{\rm )}
$r > 0.$ Moreover{\rm ,} the following quantities are equivalent
for any fixed $r > 0$ and $p \ge 1${\rm :}
\begin{equation*}
\|\hat{\mu}_{r}\|_{\infty} = \sup \left\{\frac{\mu(E(z,r))}{|
E(z,r) |^{1 - \lambda}}\hbox{\rm :}\ z \in \Omega \right\}
\end{equation*}
{\rm and}
\begin{equation*}
\| \mu\|_{p} = \sup \left\{\frac{\int_{\Omega}|f(z) |^{p} {\rm
d}\mu(z)}{\int_{\Omega}|f(z) |^{p} {\rm d}V_{\lambda}(z)}\hbox{\rm
:}\ f \in A^{p}_{\lambda}(\Omega, {\rm d}V_{\lambda}) \right\}.
\end{equation*}
\end{lem}

\begin{lem}\hskip -.3pc{\rm \cite{berger1}.}\ \
For any $r > 0,$ there exists a constant $C$ {\rm (}depending only
on $r${\rm )} such that
\begin{equation*}
C^{-1} \le \;\;|E(a,r)||k_{a}(z) |^{2}\; \le C,
\end{equation*}
for all $a \in \Omega$ and $z \in E(a,r). $

Note that if we take $z = a$ in the above estimate{\rm ,} then we
get
\begin{equation*}
C^{-1} \le \;\;|E(a,r)| K(a,a) \le C,
\end{equation*}
for all $a \in \Omega. $
\end{lem}

\begin{lem}\hskip -.3pc{\rm \cite{berger1}.}\ \
For any $r > 0, s > 0, R > 0,$ there exists $a $ constant $ C $
{\rm (}depending on $r,s,R${\rm )} such that
\begin{equation*}
C^{-1} \le \frac{|E(a,r)|}{|E(a,s)|} \le C
\end{equation*}
for all $a, b$ in $\Omega$ with $\beta(a,b) \le R. $
\end{lem}

\begin{lem}\hskip -.3pc{\rm \cite{berger2}.}\ \
For any $r >0,$ there exists a sequence $\{a_{n}\} $ in $ \Omega $
satisfying the following two conditions{\rm :}
\begin{enumerate}
\renewcommand\labelenumi{\rm (\arabic{enumi})}
\leftskip .15pc
\item $\Omega = \bigcup_{n = 1}^{\infty} E(a_{n},r); $

\item There is a positive integer $N$ such that each point $z $
in $\Omega$ belongs to at most $N$ of the sets $E(a_{n},2r).
$\vspace{-.8pc}
\end{enumerate}
\end{lem}

\begin{lem}\hskip -.3pc{\rm \cite{berger2}.}\ \
For any $r > 0$ and $p \ge 1,$ there exists a constant $ C $ {\rm
(}depending only on $r${\rm )} such that
\begin{equation*}
|f(a)|^{p} \le \frac{C}{| E(a,r) |} \int_{E(a,r)}|f(z) |^{p} {\rm
d}V(z),
\end{equation*}
for all $f$ holomorphic and $a$ in $\Omega. $
\end{lem}

\section{Bounded and compact weighted composition operators}

\setcounter{theore}{0}

\setcounter{equation}{0}

By using \cite{halmos}, we can prove the following lemma.

\begin{lem}
Take $1 \le p < \infty.$ Let $\Omega$ be a bounded symmetric
domain in $\mathbf{C}^{n},$ and $\lambda < \varepsilon_{\Omega}. $
Suppose $\varphi, \psi$ be holomorphic functions defined on
$\Omega$ such that $\varphi(\Omega) \subseteq \Omega.$ Also{\rm ,}
suppose that the weighted composition operator $W_{\varphi,\psi}$
is bounded on $ A^{p}_{\lambda}(\Omega, {\rm d}V_{\lambda})$. Then
$ W_{\varphi,\psi}$ is compact on $A^{p}_{\lambda}(\Omega, {\rm
d}V_{\lambda})$ if and only if for any norm bounded sequence
$\{f_{n}\}$ in $A^{p}_{\lambda}(\Omega, {\rm d}V_{\lambda})$ such
that $f_{n} \rightarrow 0$ uniformly on compact subsets of
$\Omega, $ then we have\break $\|W_{\varphi,\psi}
f_{n}\|_{A^{p}_{\lambda}} \rightarrow 0. $
\end{lem}

\begin{lem}
Suppose $\{f_{n}\}$ is a sequence in $A^{p}_{\lambda}(\Omega, {\rm
d}V_{\lambda})$ such that $f_{n} \rightarrow 0$ weakly in
$A^{p}_{\lambda}(\Omega, {\rm d}V_{\lambda}).$ Then $\{f_{n}\}$ is
a norm bounded sequence and $f_{n} \rightarrow 0$ uniformly on
each compact subsets of $\Omega. $
\end{lem}
Let $\varphi$ be a holomorphic function defined on $\Omega$ such
that $\varphi(\Omega) \subseteq \Omega.$ Suppose $\psi \in
A^{p}_{\lambda}(\Omega, {\rm d}V_{\lambda})$ and $\lambda <
\varepsilon_{\Omega}.$ Then the nonnegative measure
$\mu_{\varphi,\psi,p}$ is defined as
\begin{equation*}
\mu_{\varphi,\psi,p}(E) = \int_{\varphi^{-1}(E)}|\psi |^{p}{\rm
d}V_{\lambda},
\end{equation*}
where $E$ is a measurable subset of $\Omega.$

Using Lemma~2.1 of \cite{con}, we can prove the following result.

\begin{lem}
Take $1 \le p < \infty.$ Let $\Omega$ be a bounded symmetric
domain in $\mathbf{C}^{n},$ and $\lambda < \varepsilon_{\Omega}. $
Let $\varphi$ be a holomorphic mapping from $\Omega $ into
$\Omega$ and $\psi \in A^{p}_{\lambda}(\Omega, {\rm
d}V_{\lambda}).$ Then
\begin{equation*}
\int_{\Omega} g {\rm d}\mu_{\varphi,\psi,p} =
\int_{\Omega}|\psi|^{p} (g \circ \varphi) {\rm d}V_{\lambda},
\end{equation*}
where $g$ is an arbitrary measurable positive function on
$\Omega.$
\end{lem}

\begin{theor}[\!]
Take $1 \le p < \infty.$ Let $\Omega$ be a bounded symmetric
domain in $\mathbf{C}^{n},$ and $\lambda < \varepsilon_{\Omega}.$
Suppose $\varphi, \psi$ are holomorphic functions defined on
$\Omega$ such that $\varphi(\Omega) \subseteq \Omega.$ Then the
weighted composition operator $W_{\varphi,\psi}$ is bounded on
$A^{p}_{\lambda}(\Omega, {\rm d}V_{\lambda})$ if and only if the
measure $\mu_{\varphi,\psi,p}$ is a Carleson measure.
\end{theor}

\begin{proof} Suppose $W_{\varphi,\psi}$ is bounded
on $A^{p}_{\lambda}(\Omega, {\rm d}V_{\lambda}).$ Then there
exists a constant $M > 0$ such that
\begin{equation*}
\|W_{\varphi,\psi}f\|_{A^{p}_{\lambda}} \le M \|f
\|_{A^{p}_{\lambda}},
\end{equation*}
for all $f \in A^{p}_{\lambda}(\Omega, {\rm d}V_{\lambda}).$

Take
\begin{equation*}
g(\omega) = \;\;|k_{z}(\omega) |^{\frac{2 (1 - \lambda)}{p}}.
\end{equation*}
From Lemma~2.1 and using the fact that $K(\omega,z) \neq 0 $ for
any $z$ and $\omega$ in $\Omega,$ we get that $ g \in
A^{p}_{\lambda}(\Omega, {\rm d}V_{\lambda}). $ Thus, we have
\begin{align*}
\int_{\Omega}|k_{z}(\omega) |^{2 (1 - \lambda)}
{\rm d}\mu_{\varphi,\psi,p}(\omega) &= \int_{\Omega}|g(\omega) |^{p} {\rm d}\mu_{\varphi,\psi,p}(\omega)\\[.3pc]
&\le \int_{\Omega}|\psi(\omega) |^{p}|g(\varphi(\omega)) |^{p} {\rm d}V_{\lambda}(\omega)\\[.3pc]
&= \|W_{\varphi,\psi}g\|_{A^{p}_{\lambda}}^{p} < \infty.
\end{align*}
By Lemma~2.7, there exists a constant $C > 0$ (depending only on
$r$) such that
\begin{align*}
\frac{\mu_{\varphi,\psi,p}(E(z,r))}{|E(z,r)|^{1 - \lambda}} &=
\frac{1}{| E(z,r) |^{1 - \lambda}} \int_{E(z,r)} {\rm d}\mu_{\varphi,\psi,p}(\omega)\\[.4pc]
&\le C \;\;\int_{E(z,r)}|k_{z}(\omega) |^{2(1 - \lambda)}{\rm
d}\mu_{\varphi,\psi,p}(\omega) < \infty,
\end{align*}
for any $z \in \Omega.$

Thus, by using Lemma~2.6, we get that $\mu_{\varphi,\psi,p} $ is a
Carleson measure on $A^{p}_{\lambda}(\Omega, {\rm d}V_{\lambda}).$

Conversely, suppose $\mu_{\varphi,\psi,p}$ is a Carleson measure
on $A^{p}_{\lambda}(\Omega, {\rm d}V_{\lambda}).$ Then by
Theorem~1.1, there exists a constant $M > 0$ such that
\begin{equation*}
\int_{\Omega}|f(z) |^{p} {\rm d}\mu_{\varphi,\psi,p}(z) \le
\;\;M^{p} \|f\|_{A^{p}_{\lambda}}^{p},
\end{equation*}
for all $f \in A^{p}_{\lambda}(\Omega, {\rm d}V_{\lambda}). $

Hence, using Lemma~3.3, we have
\begin{align*}
\|W_{\varphi,\psi}f\|_{A^{p}_{\lambda}} &= \int_{\Omega}|\psi(z) |^{p}|f(\varphi(z)) |^{p} {\rm d}V_{\lambda}(z)\\[.4pc]
&= \int_{\Omega}| f(\omega) |^{p} {\rm d}\mu_{\varphi,\psi,p}(\omega)\\[.4pc]
&\le M^{p} \int_{\Omega}| f(\omega)|^{p} {\rm
d}V_{\lambda}(\omega) < \infty,
\end{align*}
for all $f \in A^{p}_{\lambda}(\Omega, {\rm d}V_{\lambda}).$

Hence, $W_{\varphi,\psi}$ is a bounded operator on
$A^{p}_{\lambda}(\Omega, {\rm d}V_{\lambda}).$
\end{proof}

\begin{theor}[\!]
Take $1 \le p < \infty.$ Let $\Omega$ be a bounded symmetric
domain in $\mathbf{C}^{n},$ and $\lambda < \varepsilon_{\Omega}.$
Suppose $\varphi, \psi$ are holomorphic functions defined on
$\Omega$ such that $\varphi(\Omega) \subseteq \Omega.$ Then the
weighted composition operator $W_{\varphi,\psi}$ is compact on
$A^{p}_{\lambda}(\Omega, {\rm d}V_{\lambda})$ if and only if the
measure $\mu_{\varphi,\psi,p}$ is a vanishing Carleson measure.
\end{theor}

\begin{proof} Suppose $W_{\varphi,\psi}$ is compact
on $A^{p}_{\lambda}(\Omega, {\rm d}V_{\lambda}).$ Since, $k_{z}^{1
- \lambda} \rightarrow 0$ weakly in $A^{p}_{\lambda}(\Omega, {\rm
d}V_{\lambda})$ as $z \rightarrow
\partial \Omega$, by Lemmas~3.1 and 3.2, we have
\begin{equation}
\| W_{\varphi,\psi}k_{z}^{1 - \lambda}\|_{A^{p}_{\lambda}}
\rightarrow 0 \;\; {\rm as} \;\; z \rightarrow \partial \Omega.
\end{equation}
Again, by using Lemma~2.7, we get a constant $C > 0 $ (depending
only on $r$) such that
\begin{align*}
\frac{\mu_{\varphi,\psi,p}(E(z,r))}{| E(z,r) |^{1 - \lambda}} &= \frac{1}{| E(z,r) |^{1 - \lambda}} \int_{E(z,r)} {\rm d}\mu_{\varphi,\psi,p}(\omega)\\[.4pc]
&\le C \int_{E(z,r)}|k_{z}(\omega) |^{2(1 - \lambda)}{\rm d}\mu_{\varphi,\psi,p}(\omega)\\[.4pc]
&\le C \int_{E(z,r)}|\psi(\omega) |^{2}|k_{z}(\varphi(\omega)) |^{2(1 - \lambda)}{\rm d}\mu_{\varphi,\psi,p}(\omega) \\[.4pc]
&\le C \|W_{\varphi,\psi}k_{z}^{1 - \lambda}\|_{A^{2}_{\lambda}}.
\end{align*}
Thus, condition (3.1) implies that $\mu_{\varphi,\psi,p} $ is a
vanishing Carleson measure on $A^{p}_{\lambda}(\Omega, {\rm
d}V_{\lambda}). $

Conversely, suppose that $\mu_{\varphi,\psi,p}$ is a vanishing
Carleson measure on $A^{p}_{\lambda}(\Omega, {\rm d}V_{\lambda}).$
Then, by Theorem~1.2, we have
\begin{equation}
\lim_{z \rightarrow \partial \Omega}
\frac{\mu_{\varphi,\psi,p}(E(z,r))}{| E(z,r) |^{1 - \lambda}} = 0
\end{equation}
for all (or some) $r > 0. $

So, from Theorem~3.4, we conclude that $W_{\varphi,\psi}$ is
bounded on $A^{p}_{\lambda}(\Omega, {\rm d}V_{\lambda}).$ We will
prove that $W_{\varphi,\psi}$ is a compact operator.

Let $\{f_{n}\}$ be a sequence in $ A^{p}_{\lambda}(\Omega, {\rm
d}V_{\lambda})$ and $f_{n} \rightarrow 0$ weakly. By Lemma~3.2,
$f_{n}$ is a norm bounded sequence and $f_{n}\rightarrow 0$
uniformly on each compact subsets of $\Omega.$

By Lemmas~2.7 and 2.10, there exists a constant $C > 0 $
(depending only on $r$) such that
\begin{align*}
|f(z)|^{p} &\le \frac{C}{| E(z,r) |} \int_{E(z,r) }|f(\omega)|^{p}
{\rm d}V(\omega)\\[.4pc]
&= \frac{C C_{\lambda}}{| E(z,r)|} \int_{E(z,r)}|f(\omega) |^{p}
\frac{{\rm
d}V_{\lambda}(\omega)}{K(\omega,\omega)^{\lambda}}\\[.4pc]
&\le \frac{C'}{| E(z,r) |} \int_{E(z,r)}|f(\omega) |^{p}
\frac{{\rm d}V_{\lambda}(\omega)}{| E(\omega,r) |^{-
\lambda}}\\[.4pc]
&\le C'' \frac{1}{| E(z,r) |^{1 - \lambda }}
\int_{E(z,r)}|f(\omega) |^{p} {\rm d}V_{\lambda}(\omega),
\end{align*}
for all $z \in \Omega,$ where the last inequality is obtained by
using Lemma~2.8. By combining the above inequality with Lemma~2.8,
we get
\begin{align*}
\hskip -4pc \sup \{| f(z) |^{p}{\rm :}\ z \in E(a,r)\} &\le \sup
\left\{\frac{C_{1}}{|E(z,r)|^{1 - \lambda}}
\int_{E(a,2r)}|f(\omega)
|^{p} {\rm d}V_{\lambda}(\omega){\rm :}\ z \in E(a,r)\right\}\\[.4pc]
\hskip -4pc &\le \frac{C_{2}}{|E(a,r)|^{1 - \lambda}}
\int_{E(a,2r)}|f(\omega)|^{p} {\rm d}V_{\lambda}(\omega),
\end{align*}
where $C_{1}$ and $C_{2}$ are constants depending only on
$\lambda$ and $r$. Hence, we have
\begin{align*}
\| W_{\varphi,\psi}f_{n}\|_{A^{p}_{\lambda}} &=
\int_{\Omega}|\psi(z) |^{p}|f_{n}(\varphi(z)) |^{p}{\rm d}V_{\lambda}(\omega) \\[.4pc]
&= \int_{\Omega} |f_{n}(\omega) |^{p} {\rm d}\mu_{\varphi,\psi,p} \\[.4pc]
&\le \sum_{i = 1}^{\infty} \int_{E(a_{i},r)} |f_{n}(\omega) |^{p} {\rm d}\mu_{\varphi,\psi,p} \\[.4pc]
&\le \sum_{i = 1}^{\infty} \mu_{\varphi,\psi,p}((E(a_{i},r)) \sup \{| f_{n}(\omega) |^{p}{\rm :}\ \omega \in E(a_{i},r)\} \\[.4pc]
&\le C \; \sum_{i =
1}^{\infty}\frac{\mu_{\varphi,\psi,p}((E(a_{i},r))}{| E(a_{i},r)
|^{1 - \lambda}} \int_{E(a_{i},2r)}|f_{n}(\omega) |^{p} {\rm
d}V_{\lambda}(\omega),
\end{align*}
where the sequence $\{a_{n}\}$ is the same as in Lemma~2.9.

From condition (3.2), for every $\epsilon > 0,$ there exists a
$\delta > 0 \;\;(0 < \delta < 1)$ such that
\begin{equation*}
\frac{\mu_{\varphi,\psi,p}(E(z,r))}{| E(z,r) |^{1 - \lambda}} <
\epsilon,
\end{equation*}
whenever $d(z, \partial \Omega) < \delta.$

Take
\begin{equation*}
\Omega_{1} = \{\omega \in \Omega{\rm :}\ d(z, \partial \Omega) \ge
\delta,\;\; \beta(z,\omega) \le r\}.
\end{equation*}
Then, $\Omega_{1}$ is a compact subset of $\Omega$ (as $0 < 2r <
\frac{\delta}{2}).$ Thus, we have
\begin{align*}
\hskip -4pc \int_{\Omega}|W_{\varphi,\psi}f_{n}(z) |^{p} {\rm
d}V_{\lambda}(z) &= \int_{\Omega}|f_{n}(\omega) |^{p} {\rm d}\mu_{\varphi,\psi,p}(\omega) \\[.4pc]
\hskip -4pc &\le C \; \sum_{i = 1}^{\infty}\frac{\mu_{\varphi,\psi,p}((E(a_{i},r))}{| E(a_{i},r) |^{1 - \lambda}} \int_{E(a_{i},2r)}|f_{n}(\omega) |^{p} {\rm d}V_{\lambda}(\omega)\\[.4pc]
\hskip -4pc &\le C \; \sum_{d(a, \partial \Omega) \ge \delta} \frac{\mu_{\varphi,\psi,p}((E(a_{i},r))}{| E(a_{i},r) |^{1 - \lambda}} \int_{E(a_{i},2r)}|f_{n}(\omega) |^{p} {\rm d}V_{\lambda}(\omega)\\[.4pc]
\hskip -4pc &\quad\, + C \; \sum_{d(a, \partial \Omega) < \delta}\frac{\mu_{\varphi,\psi,p}((E(a_{i},r))}{| E(a_{i},r) |^{1 - \lambda}} \int_{E(a_{i},2r)}|f_{n}(\omega) |^{p} {\rm d}V_{\lambda}(\omega)
\end{align*}
\begin{align*}
\hskip -4pc \phantom{\int_{\Omega}|W_{\varphi,\psi}f_{n}(z) |^{p}
{\rm d}V_{\lambda}(z)} &\le \textit{CN} \; \sum_{d(a, \partial \Omega) \ge \delta} \int_{E(a_{i},2r)}|f_{n}(\omega) |^{p} {\rm d}V_{\lambda}(\omega)\\[.4pc]
\hskip -4pc &\quad\, + C  \epsilon \; \sum_{d(a, \partial \Omega) < \delta} \int_{E(a_{i},2r)}|f_{n}(\omega) |^{p} {\rm d}V_{\lambda}(\omega)\\[.4pc]
\hskip -4pc &\le \textit{CNK} \; \int_{\Omega}|f_{n}(\omega) |^{p}
{\rm d}V_{\lambda}(\omega) + C \epsilon K \;
\int_{\Omega}|f_{n}(\omega) |^{p} {\rm d}V_{\lambda}(\omega),
\end{align*}
where the constant $N$ is the same as in Lemma 2.9. Since, $f_{n}$
is a norm bounded sequence and $f_{n}$ converges to zero uniformly
on compact subsets of $\Omega,$ we have $\|
W_{\varphi,\psi}f_{n}\|_{A^{p}_{\lambda}} \rightarrow 0$ as $n
\rightarrow \infty.$

The following theorem follows by combining Theorems~1.1, 1.2, 3.4
and 3.5.
\end{proof}

\begin{theor}[\!]
Take $1 \le p < \infty$ and let $\lambda < \varepsilon_{\Omega}.$
Suppose $\varphi, \psi$ be holomorphic functions defined on
$\Omega$ such that $\varphi(\Omega) \subseteq \Omega.$
Then\vspace{-.3pc}
\begin{enumerate}
\renewcommand\labelenumi{\rm (\arabic{enumi})}
\leftskip .15pc
\item The weighted composition operator $W_{\varphi,\psi}$ is
bounded on $A^{p}_{\lambda}(\Omega, {\rm d}V_{\lambda})$ if and
only if
\begin{equation*}
\hskip -1.25pc \sup_{z \in \Omega} \int_{\Omega}|k_{z}(\omega)
|^{2(1 - \lambda)} {\rm d}\mu_{\varphi,\psi,p}(\omega) < \infty.
\end{equation*}
\item The weighted composition operator $W_{\varphi,\psi}$ is
compact on $A^{p}_{\lambda}(\Omega, {\rm d}V_{\lambda})$ if and
only if
\begin{equation*}
\hskip -1.25pc \lim_{z \in \partial \Omega}
\int_{\Omega}|k_{z}(\omega) |^{2(1 - \lambda)} {\rm
d}\mu_{\varphi,\psi,p}(\omega) = 0.
\end{equation*}
\end{enumerate}
\end{theor}

\section{Schatten class weighted composition operators}

\setcounter{theore}{0}

\setcounter{equation}{0}

\begin{theor}[\!]
Take $1 \le p < \infty$ and let $\lambda < \varepsilon_{\Omega}.$
Suppose $\varphi, \psi$ are holomorphic functions defined on
$\Omega$ such that $\varphi(\Omega) \subseteq \Omega.$ Let
$W_{\varphi,\psi}$ be compact on $A^{2}_{\lambda}(\Omega, {\rm
d}V_{\lambda}).$ Then $W_{\varphi,\psi} \in S_{p}$ if and only if
the Berezin symbol of measure $\mu_{\varphi,\psi,p}$ belongs to
$L^{{p}/{2}}(\Omega, {\rm d}\sigma),$ where ${\rm d}\sigma(z) =
K(z,z) {\rm d}V(z).$
\end{theor}

\begin{proof} For any $f,g \in A^{2}_{\lambda}(\Omega, {\rm
d}V_{\lambda}),$ we have
\begin{align*}
\langle (W_{\varphi,\psi})^{*}(W_{\varphi,\psi}) f,\; g
\rangle_{A^{2}_{\lambda}} &= \langle (W_{\varphi,\psi})f, \; (W_{\varphi,\psi})g \rangle_{A^{2}_{\lambda}}\\[.4pc]
&=\int_{\Omega} f(\varphi(z)) \overline{g(\varphi(z))}|\psi(z) |^{2} {\rm d}V_{\lambda}(z) \\[.4pc]
&=\int_{\Omega} f(\omega) \overline{g(\omega)}|\psi(z) |^{2} {\rm
d}\mu_{\varphi,\psi,p}(\omega).
\end{align*}
Consider the Toeplitz operator
\begin{equation*}
T_{\mu_{\varphi,\psi,p}}f(z) = \int_{\Omega}
f(\omega)K_{\lambda}(z,\omega) \; {\rm
d}\mu_{\varphi,\psi,p}(\omega).
\end{equation*}
Then, by using Fubini's theorem, we have
\begin{align*}
\langle T_{\mu_{\varphi,\psi,p}}f, g \rangle_{A^{2}_{\lambda}} &=
\int_{\Omega} \int_{\Omega} f(\omega) K_{\lambda}(z,\omega) {\rm d}\mu_{\varphi,\psi,p}(\omega)\overline{g(z)} \; {\rm d}V_{\lambda}(z) \\[.4pc]
&=\int_{\Omega} f(\omega) \overline{\int_{\Omega} g(z) K_{\lambda}(\omega,z) {\rm d}V_{\lambda}(z)}\; {\rm d}\mu_{\varphi,\psi,p}(\omega) \\[.4pc]
&=\int_{\Omega} f(\omega) \overline{g(\omega)} \; {\rm
d}\mu_{\varphi,\psi,p}(\omega).
\end{align*}
Thus, $T_{\mu_{\varphi,\psi,p}} =
(W_{\varphi,\psi})^{*}(W_{\varphi,\psi}).$ Also, by definition an
operator $T \in S_{p}$ if and only if $(T^{*}T)^{{p}/{2}}$ is in
the trace class, and this is equivalent to saying that $T^{*} T
\in S_{{p}/{2}}$ (Lemma~1.4.6, p.~18 of \cite{zhu}). Thus
$W_{\varphi,\psi} \in S_{p}$ if and only if
$T_{\mu_{\varphi,\psi,p}} \in S_{{p}/{2}}.$ Again, by Theorem~C of
\cite{zhu1}, $T_{\mu_{\varphi,\psi,p}} \in S_{{p}/{2}}$ if and
only if the Berezin symbol of the measure $\mu_{\varphi,\psi,p}$
belongs to $L^{{p}/{2}}(\Omega, \hbox{d}\sigma).$
\end{proof}

\begin{theor}[\!]
Take $2 \le p < \infty$ and let $\lambda < \varepsilon_{\Omega}.$
Take $W_{\varphi,\psi}$ the weighted composition operator on
$A^{2}_{\lambda}(\Omega, {\rm d}V_{\lambda}).$ If
$W_{\varphi,\psi} \in S_{p},$ then
\begin{equation*}
\int_{\Omega}|\psi(z)|\left[\frac{K(\varphi(z),\varphi(z))}{K(z,z)}
\right]^{\frac{p(1 - \lambda)}{2}} {\rm d}\sigma(z) < \infty,
\end{equation*}
where ${\rm d} \sigma(z) = K(z,z) {\rm d}V(z).$
\end{theor}

\begin{proof}
By using Lemmas 2.2--2.5, we get
\begin{align*}
{\rm tr}(| W_{\varphi,\psi}|^{p}) & = {\rm tr}(| W_{\varphi,\psi}^{*}|^{p}) = {\rm tr}(| W_{\varphi,\psi}W_{\varphi,\psi}^{*}|)^{{p}/{2}} \\[.4pc]
&= \int_{\Omega} \langle (W_{\varphi,\psi}W_{\varphi,\psi}^{*})^{{p}/{2}} k_{z}^{1 - \lambda}, k_{z}^{1 - \lambda} \rangle_{A^{2}_{\lambda}} {\rm d}\sigma(z) \\[.4pc]
&\ge \int_{\Omega} \langle (W_{\varphi,\psi}W_{\varphi,\psi}^{*}) k_{z}^{1 - \lambda}, k_{z}^{1 - \lambda} \rangle_{A^{2}_{\lambda}}^{{p}/{2}} {\rm d}\sigma(z) \\[.4pc]
&= \int_{\Omega}\|W_{\varphi,\psi}^{*} k_{z}^{1 - \lambda}\|_{A^{2}_{\lambda}}^{{p}} {\rm d}\sigma(z)\\[.4pc]
&=\int_{\Omega}|\psi(z)|\left[\frac{K(\varphi(z),\varphi(z))}{K(z,z)}
\right]^{\frac{p(1 - \lambda)}{2}} {\rm d}\sigma(z).
\end{align*}
\end{proof}

\begin{theor}[\!]
Take $0 < p < 2$ and let $\lambda < \varepsilon_{\Omega}.$ Take
$W_{\varphi,\psi}$ the weighted composition operator
$W_{\varphi,\psi}$ on $A^{2}_{\lambda}(\Omega, {\rm
d}V_{\lambda}).$ If
\begin{equation*}
\int_{\Omega}|\psi(z)|\left[\frac{K(\varphi(z),\varphi(z))}{K(z,z)}
\right]^{\frac{p(1 - \lambda)}{2}} {\rm d}\sigma(z) < \infty,
\end{equation*}
where ${\rm d}\sigma(z) = K(z,z) {\rm d}V(z),$ then
$W_{\varphi,\psi} \in S_{p}.$
\end{theor}

\begin{proof} By using Lemmas~2.2--2.5, we get
\begin{align*}
{\rm tr}(| W_{\varphi,\psi}|^{p}) &= \int_{\Omega} \langle (W_{\varphi,\psi}W_{\varphi,\psi}^{*})^{{p}/{2}} k_{z}^{1 - \lambda}, k_{z}^{1 - \lambda} \rangle_{A^{2}_{\lambda}} {\rm d}\sigma(z) \\[.6pc]
&\ge \int_{\Omega} \langle (W_{\varphi,\psi}W_{\varphi,\psi}^{*}) k_{z}^{1 - \lambda}, k_{z}^{1 - \lambda} \rangle_{A^{2}_{\lambda}}^{{p}/{2}} {\rm d}\sigma(z) \\[.6pc]
&= \int_{\Omega}\|W_{\varphi,\psi}^{*} k_{z}^{1 - \lambda}\|_{A^{2}_{\lambda}}^{{p}} {\rm d}\sigma(z)\\[.6pc]
&=
\int_{\Omega}|\psi(z)|\left[\frac{K(\varphi(z),\varphi(z))}{K(z,z)}
\right]^{\frac{p(1 - \lambda)}{2}} {\rm d}\sigma(z) < \infty.
\end{align*}
Thus, $W_{\varphi,\psi}\in S_{p}.$
\end{proof}

\end{document}